\documentclass[12pt]{amsart}
\usepackage{amsmath,amssymb,amsthm,latexsym,graphicx,cite}
\def\CC{\Bbb{C}}
\def\RR{\Bbb{R}}
\def\SS{\Bbb{S}}

\newtheorem{theorem}{Theorem}

\newtheorem{corollary}[theorem]{Corollary}

\newtheorem{lemma}[theorem]{Lemma}

\newtheorem{remark}[theorem]{Remark}

\begin{document}
\title[Identities for $sin(x)$]{Identities for $\sin x$ that came from medical imaging}
\author{Peter Kuchment}
\address{Mathematics Department \\ Texas A\& M University\\ College Station, TX 77843-3368}
\email{kuchment@math.tamu.edu}
\urladdr{http://www.math.tamu.edu/\~{}kuchment}
\thanks{Research of P.~K. was supported in part by the NSF DMS Grant 0604778 and IAMCS}
\author{Sergey Lvin}
\address{
Department of Mathematics and Statistics\\ University of Maine\\ Orono, ME 04469-5257}
\email{Lvin@math.umaine.edu}
\dedicatory{Dedicated to the memory of Professor Leon Ehrenpreis,\\ a great mathematician and human being.}
\date{}
\begin{abstract}
The article describes interesting nonlinear differential identities satisfied by standard exponential
and trigonometric functions, which appeared as byproducts of medical imaging research. They look like some kind of non-commutative binomial formulas. A brief description of the origin of these identities is provided, as well as their direct algebraic derivation. Relations with separate analyticity theorems in several complex variables and some open problems are also mentioned.
\end{abstract}
\maketitle

\section{Introduction}

Working quite some time ago on some problems of medical imaging \cite{KuLv1,KuLv2}, the authors were surprised when their work produced a series of seemingly new (at least to them) identities for functions about which ``everything is known,'' such as standard exponentials and trigonometric and hyperbolic sines and cosines. It was a little bit depressing when our brave attempt to prove these identities by elementary means, without using the mathematical techniques common in tomography (c'mon, what can be difficult about $\sin x$, after all?), stalled for a while. We have succeeded eventually, albeit as we hope to persuade the reader, there are still things that are not clear about this whole business. Along the way, while trying to understand the identities, an interesting theorem of complex analysis in several variables was discovered in the joint work \cite{AEK} of one of the authors and his student with L. Ehrenpreis. We dedicate this article to the memory of Professor Leon Ehrenpreis, who passed away in August 2010 (see \cite{Notices,Struppa}).

In the next Section \ref{S:imaging} we describe briefly the medical imaging origin of the results we will be discussing. Section \ref{S:result} contains the formulation of the identities, which are then proved in the next three Sections. The paper ends with the remarks section, where in particular some open questions and relations to several complex variables theory are described.

%%%%%%%%%%%%%%%%%%%%%%%%%%%%%%%%%
\section{A brief story about medical imaging, Radon transform, and all that}\label{S:imaging}
%%%%%%%%%%%%%%%%%%%%%%%%%%%%%%%%%%%%

In tomographic medical imaging one sends through the patient's
body some kind of radiation/waves ($X$-rays, ultrasound waves, electromagnetic waves, etc.), measures the transmitted and/or reflected signals, and tries to recover from this data the internal distribution of some important for medical diagnostics internal parameter of the body (e.g., density, electrical conductivity, blood oxygenation, etc.).

The reader must have heard of at least some tomographic procedures (the best-known ones are the $X$-ray CT scan and MRI). However, not everyone might realize the heavy involvement of sophisticated, beautiful, and challenging
mathematics in creating tomographic images. In spite of its 40 years history, computed tomography still produces a steady stream of exciting mathematical problems (see, e.g., \cite{Banff,Leon,Epstein,Fokas,Feeman,GelfVil,GGG,Herman,Helg,Helg_new,Kak,KuchRadon,KuchQuinto,biomed,Natt4,Natt2001,OlafQuinto,Pal_book,Tuchin,IOut1,IOut2,vodinh,WangPATbook,WangTextbook} for description and mathematics of established and newly developing imaging modalities).

In the oldest and best established modality, X-ray CT, one deals with the problem of recovering a function $f(x)$ that is, roughly, the tissue density at the internal point $x$. Since in many cases different tissues have different densities, this function carries medically important information about the interior of the patient's body. If $f(x)$ is found, its density plot provides a picture, \textbf{tomogram} that shows to a doctor how an internal slice of the patient's body ``looks like.'' Mathematically speaking, the data obtained by the X-ray CT scanner provides the practitioner with line integrals of the sought for function $f$ (see, e.g., \cite{Epstein,Feeman,Kak,Natt4,biomed} for detailed explanation). In other words, the data one measures is the so-called \textbf{$2D$ Radon transform} that takes function $f(x)$ on the plane\footnote{We will not be precise about
conditions on functions we deal with here. It is safe to assume that they are ``very
nice'': differentiable as many times as you want and vanishing at large distances.} and integrates it over all lines $L$:
$$
f(x)\mapsto Rf(L)=\int\limits_L f(x)dl.
$$
We will parametrize lines by their \textbf{normal coordinates} $(\omega,s)$, where $\omega$ is a unit vector normal to $L$ and $s$ is the (signed) distance from the origin to $L$. Then the \textbf{normal equation} of the line $L$ has the form
$$
x\cdot \omega =s,
$$
where $x\cdot \omega$ is the
standard inner product. Thus, the Radon transform $Rf$ of a function
$f$ can be written as follows:
\begin{equation}\label{E:radon}
Rf(\omega,s)=g(\omega,s):=\int\limits_{x\cdot \omega=s}f(x)dl=
\int\limits_{-\infty}^{\infty}f(s\omega+t\omega^\perp)dt.
\end{equation}
Here $\omega^\perp=(\omega_2,-\omega_1)$ is the rotated $90^o$ vector $\omega$.

The Radon transform, besides being extremely useful in tomography, arises in a large number of other areas, e.g. in PDEs, group representation theory, etc., and thus has been the topic of thorough mathematical study \cite{Leon,GelfVil,GGG,Helg,Helg_new,Natt4,Natt2001,OlafQuinto,Pal_book}.

We are interested here in one  feature particular for the Radon transform: its range as of a linear operator between natural function spaces is rather small, i.e. of infinite co-dimension.  In other words, there are infinitely many conditions a function $g(\omega,s)$ must satisfy to be the Radon transform of a function $f$. Let us discuss them briefly.

One such condition is \textbf{evenness}:
\begin{equation}
g(-\omega,-s)=g(\omega,s)\mbox{ for all }\omega\in\SS, s\in\RR,
 \end{equation}
 where $\SS$ is the unit circle.
This reflects the fact that equations $x\cdot\omega=s$ and
$x\cdot(-\omega)=-s$ describe the same line.

More interesting are the so-called
moment conditions. To formulate them, we assume that $f$ has a compact support (in fact, decay faster than any power of $|x|$ at infinity suffices). Then for any integer $k\geq 0$ one can define \textbf{the $k$-th moment of
the data} $g(\omega,s)$:
\begin{equation}\label{E:moment}
G_k(\omega):=\int\limits_{-\infty}^\infty s^k
g(\omega,s)ds,
\end{equation}
which is a function on the circle $\SS$.

Now the \textbf{moment conditions} say that

\begin{equation}
\begin{array}{l}
\mbox{For any integer } k\geq 0,
\mbox{ the moment } G_k(\omega)\\
\mbox{can be extended from } \SS
\mbox{ to the whole plane}\\
\mbox{as a homogeneous polynomial of degree } k.
\end{array}
\end{equation}
Although it might not be immediately obvious how one discovers the moment conditions, as soon as they are formulated, it is straightforward to check their necessity.
Indeed, substituting into
(\ref{E:moment}) the definition of the Radon transform $g=Rf$, one
gets
$$
G_k(\omega):=\int\limits_{-\infty}^\infty s^k
g(\omega,s)ds=\int\limits_{-\infty}^\infty s^k \int\limits_{x\cdot
\omega=s}f(x)dlds=\int\limits_{\RR^2} (x\cdot \omega)^kf(x)dx.
$$
The last expression is clearly a homogeneous polynomial of degree
$k$ of $\omega$ (since $x\cdot \omega$ is a homogeneous linear
function of $\omega$) and thus we have established these so-called
``range conditions''.

Proving that there are no other range
conditions besides evenness and moment ones, is a different story, which
takes much more work (see, e.g. \cite{Leon,GGG,Helg,Natt4,Natt2001}).

Why would one be interested in these conditions? They are the mandatory
relations that ideal data collected from a tomographic device must
satisfy. Well, the measured data are never ideal, and so they
deviate from these conditions. Thus, knowing the range conditions
might be helpful in detecting and correcting some measurement
errors. They are of help in
other circumstances as well, for instance in completing some
missing data (e.g., see \cite{Lvin,Natt4,Natt2001,Po2} and references
therein for details). And they surely play an important role in most of
analysis of the Radon transform as an operator.

Quite a few years ago, the authors worked \cite{KuLv1,KuLv2}
on finding the range conditions for some special (weighted) Radon-type
transform arising in another popular medical imaging method, so-called
SPECT (Single Photon Emission Computed Tomography)
\cite{Fokas,Herman,KuchRadon,KuchQuinto,Natt4}. This is the so-called \textbf{exponential Radon transform}, which is defined as follows:
\begin{equation}\label{E:radon_exp}
R_\mu f(\omega,s)=g(\omega,s):=
\int\limits_{-\infty}^{\infty}f(s\omega+t\omega^\perp)e^{\mu t}dt.
\end{equation}
Here $\mu>0$ is the \textbf{attenuation coefficient}.

When we found a set of range conditions\footnote{We do not provide these conditions here, since they are
rather technical and do not make what follows easier. The interested reader can find them in \cite{KuLv1,KuLv2,KuchRadon}.} for $R_\mu$ and then proved that we had not missed any, we thought that
checking their necessity (as long as we already knew them) should be a piece of cake, just like in the example of the usual Radon transform above. You plug your transform $R_\mu f$ into the conditions and should be able to ``see'' immediately that they are satisfied. Well, when we did this, we saw an infinite series of identities for $\sin x$ (yes, the usual sine!) that we could not recognize and did not know why they should have been true. Here they are:

For any odd natural number $n$, the following identity holds:
\begin{equation}\label{E:sine_ident}
\begin{array}{c}
    \sum\limits_{k=0}^n \left(
                          \begin{array}{c}
                            n \\
                            k \\
                          \end{array}
                        \right) (\frac{d}{dx}\, -\sin x)\circ (\frac{d}{dx}\, -\sin x+i)\circ\cdots \\ \circ(\frac{d}{dx}\, -\sin x-(k-1)i)(\sin x)^{n-k}\equiv0.
\end{array}
\end{equation}
Here $i$ is the imaginary unit, $\left(\begin{array}{c}
  n \\
  k \\
\end{array}\right)$ is the binomial coefficient
``$n$ choose $k$,'' the expressions in parentheses are considered as differential operators, and $\circ$ denotes their composition. For instance,
$$
\left(\frac{d}{dx}\, -\sin x\right)u(x)=\frac{du}{dx}\, -\sin x \,\, u(x).
$$
Although we have provided the necessary proofs and moved on to doing other things, the true meaning of these identities has kept us puzzling for all these years.

In this paper we present the formulation and elementary algebraic derivation
of these identities (not only for $\sin x$, but even for some more
``elementary'' functions such as linear and exponential ones).

We would be happy if someone could shed some more light onto the
meaning and possible generalizations of these identities.

%%%%%%%%%%%%%%%%%%%%%%%%%%%%%%%%%
\section{Algebraic formulation of the identities}\label{S:result}
%%%%%%%%%%%%%%%%%%%%%%%%%%%%%%%%%%%

The identities (\ref{E:sine_ident}) can be formulated and proven in a rather general algebraic setting.
Namely, let $A$ be a commutative algebra with cancelation property and with a differentiation\footnote{I.e., $D(uv)=D(u)v+uD(v)$.} $D$ over a field $K$. An element $u\in A$ can be considered an analog of the sine function if it satisfies the equation $D^2u=\lambda^2 u$ for some $\lambda\in K$. For instance, in the case when $u=\sin x$, we are talking about the algebra $C^\infty(\RR)$ over complex numbers of all smooth functions on the variable $x$, with differentiation $D=d/dx$ and $\lambda=i$.

In what follows, the reader will not loose anything thinking, instead of the general algebraic situation, of smooth functions $f(x)$ with the usual derivative $D=d/dx$. However, the general algebraic notations ($A$, $D$, $\lambda$) lead to somewhat easier to read formulas, so we will stuck with them most of the time.

We are interested in solutions of simple differential equations
\begin{equation}\label{E:1st}
    Du=\lambda u
\end{equation}
and
\begin{equation}\label{E:2nd}
    D^2u=\lambda^2 u,
\end{equation}
where $\lambda\in K$. Thus, in the example of the function algebra with usual differentiation,
solutions of equation (\ref{E:1st}) for $\lambda=0$ are constants,
while for non zero values of $\lambda$ they give us exponential
functions. Analogously, solutions of (\ref{E:2nd}) for $\lambda=0$
are linear functions, while for non zero values of $\lambda$ they
in particular give us trigonometric and hyperbolic sine and
cosine.

We now formulate the main results.

\begin{theorem}\label{T:1st}
Any solution $u$ of the first-order equation $Du=\lambda u$
satisfies for any natural $n$ the following identity:
\begin{equation}\label{E:ident_1st}
\sum\limits_{k=0}^{n}\left(%
\begin{array}{c}
  n \\
  k \\
\end{array}%
\right)(D-u)\circ(D-u+\lambda )\circ ... \circ ( D-u+(k-1)\lambda)
u^{n-k}=0.
\end{equation}
\end{theorem}

\begin{theorem}\label{T:2nd}
Any solution $u$ of the second-order equation $D^2u=\lambda^2 u$
satisfies for any {\bf odd} natural $n$ the identity
(\ref{E:ident_1st}).
\end{theorem}

The expression $u^{n-k}$ at the very right is
just the $(n-k)$th power of $u$. Expressions
like $(D-u)$ are considered as {\bf operators} on the algebra $A$. For instance, $(D-u)$
applied to an element $f\in A$ produces
$(D-u)f=Df-uf$. The little circles $\circ$ mean
{\bf composition} of these operators from right to left. For
instance,
$$
\begin{array}{c}
(D-u)\circ (D-u+\lambda) f=(\dfrac{d}{dx}-u(x))(\dfrac{df}{dx}-u(x)f(x)+\lambda
f(x))\\=\dfrac{d^2 f}{dx^2}-\dfrac{d(uf)}{dx}+\lambda
\dfrac{df}{dx} -u(x)\dfrac{df}{dx} +u^2(x)f(x)-\lambda u(x)f(x).
\end{array}
$$
So, the expression in (\ref{E:ident_1st}) directs us to raise $u$
to the power $(n-k)$, then apply to it the operator
$(D-u+(k-1)\lambda)$, then to apply to the resulting function the
operator $(D-u+(k-2)\lambda)$, etc. To  avoid confusion in the
case of $k=0$, it might be easier to think that the composition in
(\ref{E:ident_1st}) has $k$ operator factors. In other words, the
term for $k=0$ does not have any of these factors and thus is just
$u^n$.

An important thing to notice is that, as Calculus teaches us, the
operations $D$ and $u$ (differentiation and multiplication by $u$)
do not commute, i.e. the result of their composition depends on its order.
If this were not so, the identities
(\ref{E:ident_1st}) would become trivial, as we will see a little
bit later.

Since solutions of the equations (\ref{E:1st})-(\ref{E:2nd}) are
very simple functions (constants and exponential, linear,
trigonometric or hyperbolic functions), let us look at some
examples to see whether the identities (\ref{E:ident_1st}) might
be trivial, or at least look familiar.

Take (\ref{E:1st}) for $\lambda=0$, so the solution $u$ is just a constant $C$. Then
(\ref{E:ident_1st}) looks as follows:
$$
\sum\limits_{k=0}^{n}\left(%
\begin{array}{c}
  n \\
  k \\
\end{array}%
\right)(D-C)^k C^{n-k}=0.
$$
Since $D$ acting on any constant gives zero result, this identity boils down to
$$
\sum\limits_{k=0}^{n}\left(%
\begin{array}{c}
  n \\
  k \\
\end{array}%
\right)(-C)^k C^{n-k}=(-C+C)^n=0.
$$
This is indeed easy and proves the particular case of Theorem \ref{T:1st} for $\lambda=0$.

Consider now the particular case of (\ref{E:1st}) when $\lambda=1$ and $u=e^x$. Then
(\ref{E:ident_1st}) reduces to
$$
\sum\limits_{k=0}^{n}\left(%
\begin{array}{c}
  n \\
  k \\
\end{array}%
\right)\left(\frac{d}{dx}-e^x\right)\circ\left(\frac{d}{dx}-e^x+1\right)\circ ... \circ \left(\frac{d}{dx}-e^x+(k-1)\right)
e^{(n-k)x}=0.
$$
Hmmm... Does not look familiar? Neither it did to us. Well, we will not prove this particular identity, since
we are going to prove the more general Theorems \ref{T:1st} and \ref{T:2nd}.

Let us look at an example of a solution of the second-order
equation $D^2u=\lambda^2u$, say take $u=sinx$ and $\lambda=i$.
Then the identity (\ref{E:ident_1st}) for $n=3$ becomes (using
more familiar  $\frac{d}{dx}$ instead of $D$)
\begin{equation}\label{E:sine}
\begin{array}{c}
\sin^{3}x+3(\dfrac{d}{dx}-\sin x)\sin^{2}x + 3(\dfrac{d}{dx}-\sin
x)\circ(\dfrac{d}{dx}-\sin x+i)\sin x\\+ (\dfrac{d}{dx}-\sin
x)\circ(\dfrac{d}{dx}-\sin x+i)\circ(\dfrac{d}{dx}-\sin x+2i)1=0.
\end{array}
\end{equation}
Obvious, isn't it? Well, do not despair if you do not see it right
away, we don't either.

As you will see, the direct algebraic proofs of Theorems \ref{T:1st} and \ref{T:2nd} (at
least the proofs that we know) are different when $\lambda
=0$ and when  $\lambda \neq 0$. Certainly, it is not hard to derive the
identities for $\lambda=0$ from the ones for $\lambda\neq 0$ by
taking the limit when $\lambda \rightarrow 0$.
We, however, would like to see direct algebraic proofs of the identities for both
cases $\lambda =0$ and $\lambda \neq 0$. Such proofs are provided below.

%%%%%%%%%%%%%%%%%%%%%%%%%%%%%%%%
\section{Proof of Theorem \ref{T:1st}}\label{S:Th1}
%%%%%%%%%%%%%%%%%%%%%%%%%%%

We have already proven it above for $\lambda=0$. So, let us assume that $\lambda\neq 0$.

Suppose that $u$ satisfies the first-order equation $Du=\lambda u$. We
can apply a \textbf{gauge transform} to simplify our identities. Namely, we will use
the easy to check identity
\begin{equation}\label{E:commut3}
u^{-m}(D-u)u^{m}=D-u+m\lambda.
\end{equation}
We can now use (\ref{E:commut3}) in each factor of (\ref{E:ident_1st}) to rewrite
it as follows:
\begin{equation}
\left(\sum\limits_{k=0}^{n}\left(%
\begin{array}{c}
  n \\
  k \\
\end{array}%
\right)[(D-u)\circ u^{-1}]^{k}\right)u^{n}=0.
\end{equation}
Here we notice that we deal with an operator binomial
$\sum\limits_{k=0}^{n}\left(%
\begin{array}{c}
  n \\
  k \\
\end{array}%
\right)A^{k}$, where $A=(D-u)\circ u^{-1}$. Since only one
operator is involved, no non-commutativity arises, and thus the
usual binomial formula works. This reduces identity
\eqref{E:ident_1st} to
\begin{equation}
 \lbrack (D-u)\circ
u^{-1}+1]^{n}u^{n} =(D\circ u^{-1})^{n}u^{n}=0.
\end{equation}
Now an immediate calculation using $Du=\lambda u$ shows that $(D\circ u^{-1})^{n}u^{n}=0$
holds.

\begin{remark}\label{R:negative}
A careful reader might object to our calculations in this section
that involve negative powers of $u$, as well as positive ones, especially in the abstract setting of a commutative differential algebra $A$ with cancelation. However, this algebraic problem can be overcome (see \cite{KuLv1}).
\end{remark}

%%%%%%%%%%%%%%%%%%%%%%%%%%%%%
\section{Proof of Theorem \ref{T:2nd} for $\lambda =0$}\label{S:Th2_0}
%%%%%%%%%%%%%%%%%%%%%%%%%%%%

We assume now that $\lambda=0$ and thus the function  $u$ satisfies the equation
$D^2u=0$. In this case the identities (\ref{E:ident_1st}) that we intend to prove
simplify to
\begin{equation}\label{E:start}
\sum\limits_{k=0}^{n}\left(%
\begin{array}{c}
  n \\
  k \\
\end{array}%
\right)(D-u)^{k}u^{n-k}=0
\end{equation}
for any odd natural $n$.
\begin{remark}\label{R:even}
This identity fails for $n=2$.
\end{remark}

One can understand (\ref{E:start}) also as follows:
\begin{equation}\label{E:l=0}
\left(\sum\limits_{k=0}^{n}\left(%
\begin{array}{c}
  n \\
  k \\
\end{array}%
\right)(D-u)^{k}\circ u^{n-k}\right)1=0.
\end{equation}
The expression in the parentheses resembles the familiar binomial expression
$$
\sum\limits_{k=0}^{n}\left(%
\begin{array}{c}
  n \\
  k \\
\end{array}%
\right)a^{k}b^{n-k},
$$
where $a=D-u$ and $b=u$. If the operators commuted, then
the whole thing would have been equal to $(a+b)^n=(D-u+u)^n
1=D^n 1=\dfrac{d^n 1}{dx^n}=0$. Done!

Hold on, not so fast! First of all,
since $D$ and $u$ do not commute, the usual binomial formula does not apply. Besides,
it looks like we obtained the identity (\ref{E:ident_1st}) for
{\bf all} natural $n$, which is incorrect, as we saw in
Remark \ref{R:even}. So, let us take a closer look.

If $\lambda=0$, then function $u$ is linear. Let us concentrate on
the case when $u(x)=x$.

Let us first rewrite our identity (\ref{E:l=0}) for the case
when $u=x$:
\begin{equation}\label{E:l=0,u=x}
\sum\limits_{k=0}^{n}\left(%
\begin{array}{c}
  n \\
  k \\
\end{array}%
\right)\left(D-x\right)^{k} x^{n-k}=0.
\end{equation}

We will now use another popular {\bf gauge
transformation} trick that will enable us to
rewrite (\ref{E:l=0,u=x}) in a slightly different form. It is based on the following simple
identity:
\begin{equation}\label{E:gauge}
\left(D-x\right)^{k}\left(e^{\frac{x^2}{2}}f(x)\right)
=e^{\frac{x^2}{2}}D^k f(x).
\end{equation}
Thus, commutation with the function
$e^{\frac{x^2}{2}}$ kills the term $-x$ added to the
 derivative and $\left(D-x\right)^{k}$ becomes $D^k$:
 $$
 e^{-\frac{x^2}{2}}\circ \left(D-x\right)^{k}\circ e^{\frac{x^2}{2}}
 =D^k.
 $$
This implies the following
\begin{lemma}\label{L:gauge}
The identities (\ref{E:l=0,u=x}) are equivalent to
\begin{equation}\label{E:gauge_ident}
\left(\sum\limits_{k=0}^{n}\left(%
\begin{array}{c}
  n \\
  k \\
\end{array}%
\right)D^k \circ x^{n-k}\right)e^{-\frac{x^2}{2}}=0.
\end{equation}
\end{lemma}

So, now the proof of Theorem \ref{T:2nd} for $\lambda=0$ boils
down to proving (\ref{E:gauge_ident}) for odd $n$. Notice that if
$D$ and $x$ commuted, then the expression in parentheses in
(\ref{E:gauge_ident}) would become just $(D+x)^n$, and thus, due
to the easily verified equality
\begin{equation}\label{E:exponent}
\left(D+x\right)e^{-\frac{x^2}{2}}=0,
\end{equation}
we would immediately obtain the validity of
(\ref{E:gauge_ident}). However, life is not so easy and thus the
operators do not commute. (Besides, as we already know,
(\ref{E:gauge_ident}) is guaranteed for odd values of $n$ only).

Based on (\ref{E:gauge_ident}), let us introduce the differential
operators
$$
P_n(D,x)=\sum\limits_{k=0}^{n}\left(%
\begin{array}{c}
  n \\
  k \\
\end{array}%
\right)D^k \circ x^{n-k}.
$$
The punch line is in the following statement:
\begin{lemma}
The following recurrence relation holds:
\begin{equation}\label{E:recurrent}
   P_{n+2}(D,x)=P_{n+1}(D,x)\circ (D+x)+(n+1)P_{n}(D,x).
\end{equation}
\end{lemma}

{\bf Proof of the Lemma}. We will compute $P_{n+1}(D,x)\circ
(D+x)$ and show that it coincides with
$P_{n+2}(D,x)-(n+1)P_{n}(D,x)$. In order to do so, we will use the
easy to verify commutation relation
$$
D\circ x^m =x^m\circ D +mx^{m-1}
$$
and its consequence
\begin{equation}\label{E:commut2}
x^m\circ D =D\circ x^m-mx^{m-1}.
\end{equation}

Now, take a deep breath, and ...
\begin{equation}\label{E:comput}
\begin{array}{c}
P_{n+1}(D,x)\circ (D+x)=\left(\sum\limits_{k=0}^{n+1}\left(%
\begin{array}{c}
  n+1 \\
  k \\
\end{array}%
\right)D^k \circ x^{n+1-k}\right)\circ (D+x)\\
=\sum\limits_{k=0}^{n+1}\left(%
\begin{array}{c}
  n+1 \\
  k \\
\end{array}%
\right)D^k \circ x^{n+1-k}\circ D+\sum\limits_{k=0}^{n+1}\left(%
\begin{array}{c}
  n+1 \\
  k \\
\end{array}%
\right)D^k \circ x^{n+2-k}.
\end{array}
\end{equation}
Let us now use (\ref{E:commut2}) to commute $x^{n+1-k}$ with $D$ in the first sum to get
\begin{equation}\label{E:comput2}
\begin{array}{c}
    \sum\limits_{k=0}^{n+1}\left(%
\begin{array}{c}
  n+1 \\
  k \\
\end{array}%
\right)D^{k+1} \circ x^{n+1-k}+\sum\limits_{k=0}^{n+1}\left(%
\begin{array}{c}
  n+1 \\
  k \\
\end{array}%
\right)D^k \circ x^{n+2-k}\\
-\sum\limits_{k=0}^{n+1}(n+1-k)\left(%
\begin{array}{c}
  n+1 \\
  k \\
\end{array}%
\right)D^{k} \circ x^{n-k}.
\end{array}
\end{equation}
The first two sums can be rewritten as
\begin{equation}\label{E:comput3}
\begin{array}{c}
    \sum\limits_{k=1}^{n+2}\left(%
\begin{array}{c}
  (n+2)-1 \\
  k-1 \\
\end{array}%
\right)D^{k+1} \circ x^{(n+2)-k}+\sum\limits_{k=0}^{(n+2)-1}\left(%
\begin{array}{c}
  (n+2)-1 \\
  k \\
\end{array}%
\right)D^k \circ x^{(n+2)-k}\\
=\sum\limits_{k=0}^{n+2}\left(%
\begin{array}{c}
  n+2 \\
  k \\
\end{array}%
\right)D^k \circ x^{(n+2)-k}=P_{n+2}(D,x).
\end{array}
\end{equation}
We used here the standard property of the Pascal
triangle:
$$
\left(
\begin{array}{c}
  (n+2)-1 \\
  k-1 \\
\end{array}\right)+\left(%
\begin{array}{c}
  (n+2)-1 \\
  k \\
\end{array}%
\right)=\left(%
\begin{array}{c}
  n+2 \\
  k \\
\end{array}%
\right).
$$
It only remains to handle the remainder
$$
\sum\limits_{k=0}^{n+1}(n+1-k)\left(%
\begin{array}{c}
  n+1 \\
  k \\
\end{array}%
\right)D^{k} \circ x^{n-k}.
$$
This expression is equal to $(n+1)P_{n}(D,x)$, proving which
only requires the knowledge of what binomial coefficients are.
This finishes the proof of the Lemma. \qed

\begin{corollary}\label{C:divisibility}
For any odd natural $n$, the following factorization holds with some
operator $Q(D,x)$:
\begin{equation}\label{E:factor}
    P_n=Q\circ (D+x).
\end{equation}
\end{corollary}
\proof Let us prove this by induction. When $n=1$, we have
$P_1=D+x$, and the statement is obvious. Assume that we have
proven it for some odd $n$, i.e. $P_n=Q\circ (D+x)$. Then the
previous Lemma implies that $P_{n+2}=\left(P_{n+1}+(n+1)Q\right)\circ
(D+x)$. \qed

\begin{remark}
The factorization $P_n=Q\circ (D+x)$ fails for $n=2$.
\end{remark}

Now the proof of Theorem \ref{T:2nd} for $\lambda=0$ is immediate.
Indeed, for any odd $n$, the left hand side in the identity
(\ref{E:gauge_ident}) in question becomes $P_n
e^{-\frac{x^2}{2}}=Q\circ (D+x) e^{-\frac{x^2}{2}}$. Applying
(\ref{E:exponent}), we conclude that this is zero. This finishes
the proof of the Theorem for $\lambda=0$ and $u=x$.

The proof provided for $u=x$ can be easily generalized to an arbitrary linear function $u=ax+b$. Alternatively, one can derive the result for any linear function from the one for $u=x$.
\qed

%%%%%%%%%%%%%%%%%%%%%%%%%%%%%
\section{Proof of Theorem \ref{T:2nd} for $\lambda \neq 0$}
%%%%%%%%%%%%%%%%%%%%%%%%%%%

Let us now outline the steps of the proof of Theorem \ref{T:2nd}
in case $\lambda \neq 0$.

\textbf{Step 1.} If $u$ satisfies the first-order equation
$Du=\lambda u$ (and thus certainly $D^2u=\lambda^2u$ as well), the statement follows from Theorem \ref{T:1st}.

\textbf{Step 2.} Suppose $Du=-\lambda u$. Then (\ref{E:ident_1st})
is true for odd $n$, because one can easily check that
``symmetric'' terms in the sum that correspond to $k$ and $n-k$ cancel out.

\textbf{Step 3.} Let us now suppose that $u$ is any solution of $D^2u=\lambda ^{2}u$. If
our identities were linear with respect to $u$ (which they are
not), then the cases when $Du=\pm\lambda u$ discussed in the two
previous steps would suffice, since any solution could be
expanded into a sum of these two. Nonlinearity seems to destroy
this idea. However, having nothing better in mind, let us still
try. Thus, $u=v+w$, where $Dv=-\lambda v$ and $Dw=\lambda w$ (notice that the
assumption $\lambda \neq 0$ is critical for the possibility of such
decomposition)\footnote{This is just standard exponential
representation of solutions studied in ODEs (or, if you will, the
representation of $\sin x$ as combination of $e^{\pm ix}$).}.

\textbf{Step 4.} Let us use the same commutation trick (\ref{E:commut3})
$$
u^{-m}(D-u)u^{m}=D-u+m\lambda
$$
as before, but
using commutation with $w$ rather than $u$.
Then the left-hand side of (\ref{E:ident_1st}) can be
written as
\begin{equation}\label{E:binom}
\begin{array}{c}
\sum\limits_{k=0}^{n}\left(%
\begin{array}{c}
  n \\
  k \\
\end{array}%
\right)(D-v-w)\circ ...\circ ( D-v-w+(k-1)\lambda ) (v+w)^{n-k}\\
=\left[ \sum\limits_{k=0}^{n}\left(\begin{array}{c}
  n \\
  k \\
\end{array}\right)[(D-v)\circ w^{-1}-1]^{k}(vw^{-1}+1)^{n-k}\right]
w^{n}.
\end{array}
\end{equation}

\textbf{Step 5.} Let us introduce the following operator notations:
$A=(D-v)\circ w^{-1}$ and $B=vw^{-1}$. Then the last sum becomes
$$
\sum\limits \left(\begin{array}{c}
  n \\
  k \\
\end{array}\right)(A-1)^{k}(B+1)^{n-k}.
$$
If the operators $A$ and $B$ commuted, then according to the binomial
formula this would boil down to $(A+B)^n$ and thus also to
$\sum\limits \left(\begin{array}{c}
  n \\
  k \\
\end{array}\right)A^{k}B^{n-k}$.
The interesting thing is that the latter conclusion holds even without
commutativity:
\begin{lemma}\label{L:noncom}
For any two operators $A,B$ the following equality holds:
$$
\sum\limits \left(\begin{array}{c}
  n \\
  k \\
\end{array}\right)(A-1)^{k}(B+1)^{n-k}=\sum\limits \left(\begin{array}{c}
  n \\
  k \\
\end{array}\right)A^{k}B^{n-k}.
$$
\end{lemma}
Indeed, this equality does not require changing the order of the factors. Thus,
if it holds for commuting operators, then it does for non-commuting
ones as well.

How can this help us with the proof of the theorem? It allows us to drop
the terms $\mp 1$ in (\ref{E:binom}) to get
$$
\left[ \sum\limits_{k=0}^{n}\left(\begin{array}{c}
  n \\
  k \\
\end{array}\right)[(D-v)\circ w^{-1}]^{k}(vw^{-1})^{n-k}\right]
w^{n}.
$$
Now one can reverse Step 4 (undoing the commutations with powers of $w$
we have done) and rewrite (\ref{E:binom}) as
$$
\sum\limits_{k=0}^{n}\left(\begin{array}{c}
  n \\
  k \\
\end{array}\right)(D-v)\circ ...\circ ( D-v+(k-1)\lambda) v^{n-k},
$$
which is equal to zero due to the Step 2, since $v$ solves the
first-order equation $Dv=-\lambda v$. This proves Theorem
\ref{T:2nd} for $\lambda\neq 0$. \qed

%%%%%%%%%%%%%%%%%%%%%%%%%%
\section{Remarks and generalizations}
%%%%%%%%%%%%%%%%%%%%%%%%%%

\begin{itemize}

\item For the readers with a taste for generalizations,
we can mention that the identities (\ref{E:ident_1st}) hold in a much more general situation.
Namely, $u$ can be assumed to be an element of a commutative algebra with cancelation and with differentiation $D$
over a field $\Lambda$. Then $\lambda$ should be an element of the field $\Lambda$ \cite{KuLv1}.

\item Generalizations of some of the identities we discussed are available:
\begin{theorem}
If $u$ satisfies $D^2u=0$, then the identity
\begin{equation}
\sum_{k=0}^{n}\left(%
\begin{array}{c}
  n \\
  k \\
\end{array}%
\right)(D-u)^{k}\circ D^{m}u^{n-k}=0
\end{equation}
holds for all odd integer $n\geq 1$ and even $m\geq 0$.
\end{theorem}

\item Several integral geometry and tomography experts devoted their time and effort to trying to understand better the meaning of these strange range conditions. This is also what one of the authors set out to do with V.~Aguilar and L.~Ehrenpreis. Surprisingly, as the result, the identities (\ref{E:ident_1st}) were related \cite{AEK} to the so-called separate analyticity theorems (Hartogs-Bernstein theorems) in several complex variables. For instance, the following amazing theorem is
essentially equivalent to these identities
\cite{AEK,Oktem1,Oktem2}:
\begin{theorem}
Let $\Omega$ be a disk in $\RR^2$ and $f(x)$ be a function in the exterior of $\Omega$. Suppose that when restricted to any tangent line $L$ to $\partial\Omega$, the function $f|_L$, as a function of one real variable extends to an entire function on the complexification of $L$. Then $f$, as a function on $\RR^2\setminus \Omega$ extends to an entire function on $\CC^2$.
\end{theorem}

Well, this fact also did not look obvious. Analyticity of $f$ in a complex neighborhood of $\RR^2\setminus D$ follows from the old (and not well known) separate analyticity theorem by S.~Bernstein (see \cite{AkhRon}), however this theorem cannot produce statement about $f$ being an entire function. Thus, since proving the above theorem, a couple of things about it kept bothering us for several years. First of all, this is a several complex variables fact, while our proof did not look like a SCV argument at all. Is there a truly complex analysis proof? Another, related, question is whether such a theorem can be proven for a different convex body instead of a disk $\Omega$? A SCV proof was later provided in \cite{Oktem1,Oktem2}, although it was rather complicated and was not generalizable (at least, easily) to other convex curves. Leon Ehrenpreis has worked out some other examples of convex algebraic curves (unpublished), but general picture remained unclear. Finally, A.~Tumanov presented recently \cite{Tumanov} a beautiful short proof based on attachment of analytic disks (where Tumanov is a great expert), which works for any strictly convex body $\Omega$ with a mild conditions on the smoothness of its boundary.

More discussion of the relations of the range conditions with complex analysis can be
found for instance in \cite{EKP,KuchRadon,KuLv1,KuLv2}.

\item It would be interesting to find an algebraic proof of the
identities that would work simultaneously for $\lambda=0$ as well
as for $\lambda\neq 0$.

\item It is clear that the identities discussed must be related to
special function theory and group representations. It would be
interesting to understand such relations.

\item  The formulations of the Theorems \ref{T:1st} and \ref{T:2nd}
lead to the natural question: what can be said for solutions of the
equation of 3rd order $D^3u=\lambda^3u$ and higher? A natural guess would be
that the same identities hold, but only for an arithmetic sequence of numbers
$n$ with the difference equal to $3$. It is hard to compute these expressions by hand even
for small values of $n$, say for $n=4$. Ms. E. Rodriguez, a former Masters student
of P.~K., has used the Maple symbolic algebra system to check this conjecture.
The result is negative, the natural conjecture fails for equations
of third order \cite{Elaina}. So, what (if anything) happens to
solutions of higher order differential equations $D^mu=\lambda^m u$?
Are they deprived of any such identities? We do not know the answer to this question.
\end{itemize}

\section*{Acknowledgments}

The authors express their gratitude to S.~Fulling, D.~Hensley,
Z.~Sunik and G.~Tee for useful comments.

\end{document}